\input amstex
\documentstyle{amsppt}
\document
\magnification=1200
\NoBlackBoxes
\nologo

\vsize18cm


\bigskip

\centerline{\bf MATHEMATICAL KNOWLEDGE:}

\smallskip

\centerline{\bf INTERNAL, SOCIAL AND CULTURAL ASPECTS\footnotemark1}
\footnotetext{Commissioned for vol. 2 of ``Mathematics
and Culture'', ed. by C.~Bartocci and P.~Odifreddi.}

\medskip

\centerline{Yu.~I.~Manin}

\medskip

\centerline{\it Max--Planck--Institut f\"ur Mathematik, Bonn, Germany,}

\centerline{\it and Northwestern University, Evanston, USA}

\bigskip

\centerline{\bf 0. Preface}

\medskip

\hfill{\it He is onto us and how, on the one hand, we take pride 
in building an elegant world utterly divorced from the demands of reality
and, on the other, claim that our ideas underlie virtually all 
technological developments of significance.}

\smallskip

\hfill{\it D.~Mumford (from the Preface to [Ens]).}

\medskip

Pure mathematics is an immense organism
built entirely and exclusively of ideas that emerge in the minds
of mathematicians and live within these minds.

\smallskip 

If one wishes to shake off the somewhat uneasy
feeling that such a statement can provoke, there are
at least three escape routes.

\smallskip

First, one can simply identify mathematics with the contents of
mathematical manuscripts, books, papers and lectures,
with the increasingly growing net of theorems, definitions,
proofs, constructions, conjectures (should I include
software as well ?...) --  with what contemporary mathematicians present 
at the conferences, keep in the libraries
and electronic archives, take pride in, award each other
for, and occasionally bitterly dispute the origin of.
In short, mathematics is simply what mathematicians are doing,
exactly in the same way as music is what musicians are doing.

\smallskip

Second, one can argue
that mathematics is a human activity deeply rooted
in reality, and permanently returning to reality.
From  counting on one's fingers to moon--landing to Google,
we are doing mathematics in order to understand, create,
and handle things, and perhaps {\it this understanding}
is mathematics rather than intangible murmur
of accompanying abstractions. Mathematicians are thus more or less
responsible actors of human history,
like Archimedes helping to defend Syracuse (and to save 
a local tyrant), Alan Turing cryptanalyzing Marshal Rommel's 
intercepted military dispatches to Berlin, or John von Neumann
suggesting high altitude detonation as an efficient tactics of bombing.
Accepting this viewpoint, mathematicians can defend
their trade by stressing its social utility. 
In this role, a mathematician can be as morally confused as the next person,
and if I were to put on display some trade--specific particularities
of such a confusion,
I could not find anything better than the bitter irony
of [B-BH] (p. 11): ``[...] mathematics can also be an indispensable tool. Thus, when the effect of fragmentation bombs on human bodies was to be tested but
{\it humanitarian concerns prohibited testing on pigs} (italics mine. Yu.~M.), mathematical simulation was put into play.''

\smallskip

Or, third, there is a noble 
vision of the great Castle of Mathematics, towering 
somewhere in the Platonic World of Ideas, which we humbly
and devotedly discover (rather than invent).
The greatest mathematicians manage to grasp
outlines of the Grand Design, but even those to whom
only a pattern on a small kitchen tile is revealed,
can be blissfully happy. Alternatively, if one is inclined to use
a semiotic metaphor, Mathematics is a proto--text
whose existence is only postulated but which
nevertheless underlies all corrupted and fragmentary
copies we are bound to deal with. The identity of the writer
of this proto--text (or of the builder of the Castle)
is anybody's guess, but Georg Cantor
with his vision of infinity of infinities directly
inspired by God, and
Kurt G\"odel with his ``ontological proof'',
seemingly had no doubts on this matter. 

\smallskip

Various shades and mixes of these three attitudes, social positions, 
and implicated choices of the individual behavior,
color the whole discussion that follows.
The only goal of this concise Preface is to make 
the reader conscious of the intrinsic tensions in our
presentation, rather than imitate clear vision
and offer definite judgements where there are none. 
 
\smallskip

One last warning about historical references
in this exposition. There are two different modes of
reading old texts: one, to understand the times and ethnos
they were written in, another - to throw some light on the
values and prejudices of our times. In the history of mathematics,
the polar attitudes are represented by ``ethnomathematics" vs
Bourbaki style history. 

\smallskip

For the sake of this presentation, I explicitly
and consciously adopt the ``modernizing'' viewpoint.

\smallskip

{\it Acknowledgement.} Silke Wimmer--Zagier provided
some sources on the history of Chinese and Japanese
mathematics and discussed their relevance to this
project. Dmitri Manin explained me Google's strategy
of page ranking. I gratefully acknowledge their generous help.

\bigskip

\centerline{\bf I. Mathematical knowledge}

\medskip

{\bf I.1. Bird's eye view.} Sir Michael Atiyah starts his
report [At] with  the following broad outline: ``The three
great branches of mathematics are, in historical order,
Geometry, Algebra and Analysis. Geometry we owe
essentially to Greek civilization, Algebra is
of Indo--Arab origin and Analysis (or the Calculus)
was the creation of Newton and Leibniz, ushering in the modern era.''
He then explains that in the realm of physics,
these branches correspond respectively to 
the (study of) Space/Time/Continuum: ``There is little
argument about Geometry being the study of space, but it
is perhaps less obvious that Algebra is the study of time.
But, any algebraic system involves the performance of sequential
operations (addition, multiplication, etc.) and these are 
conceived as being performed one after another. 
In other words Algebra requires time for
its meaning (even if we usually only need
discrete instants of time).''

\smallskip

One can argue for an alternative viewpoint on Algebra
according to which it has most intimate relations
not with Physics but with Language.
In fact, observing the graduate emergence of place--value
notation for numbers, and later algebraic notation for variables and 
operations, one can recognize two historical stages.

\smallskip 

At the first stage, notation serves principally
to shorten and unify the symbolic representation
of a certain pool of meanings. At this stage,
a natural language could (and did) serve the same goal,
only less efficiently. Therefore one can reasonably
compare this process with the development of
a specialized sub--dialect of the natural language.
The so--called Roman numerals still in use for ornamental
purposes are fossilized remnants of this stage.
As another helpful comparison, perhaps more streamlined and
better documented, one can invoke the emergence and evolution
of chemical notation.

\smallskip

At the second stage, algorithms for addition/multiplication
and later division of numbers in a place--value
notation are devised. In a parallel development,  
variables and algebraic operations 
start to be combined into identities and equations,
and then to strings of equations obeying universal rules
of identical transformations/deductions. 
At this stage, expressions in the new (mathematical) dialect
become not so much carriers of certain meanings as
a grist for the mill of computations. It is this
shift of the meaning, from the more or less explicit 
semantics of notation to
the hidden semantics of algorithms transforming 
strings of symbols, that was the crucial chain of events
marking the birth of Algebra.

\smallskip

Nothing similar to this second stage happened to
the natural languages. To the contrary, when
in the 60s of the twentieth century large computers
made possible first experiments with algorithmic
processing of texts in English, Russian, French
(e.g. for implementing automatic translation), it became
clear how unsuitable for computer processing
natural languages were. Huge data bases for
vocabularies were indispensable. Intricate and illogical
nets of rules governed morphology, word order,
and compatibility of grammatical constructions;
worse, in different languages these rules were
capriciously contradictory. And after all efforts,
automatic translation without subsequent
editing by a human being never produced satisfactory results.

\smallskip

This property of human languages -- their
resistance to algorithmic processing -- is perhaps
the ultimate reason why only mathematics
can furnish an adequate language for physics.
It is not that we lack words for expressing
all this $E=mc^2$ and $\int e^{iS(\phi )}D\phi$ stuff  -- words
can be and are easily invented -- the point is that
we still would  not be able to do anything with these
great discoveries if we had only words for them.

\smallskip

But we cannot just skip words and deal only with formulas either.
Words in mathematical and scientific texts play three
basic roles. First, they furnish multiple bridges between the 
physical reality
and the world of mathematical abstractions. Second, they
carry value judgements, sometimes explicit,
sometimes implicit, governing our choices of particular chains
of  mathematical reasonings, in the vast tree
of ``all'' feasible but mostly empty formal deductions.
And last but not least, they allow us to communicate, teach and learn.

\smallskip

I will conclude with a penetrating comment
of Paul Samuelson regarding use of words vs mathematical
symbols in economic models (cited from [CaBa]):
``When we tackle them [the problems of economic theory] by words,
we are solving the same equations as when we write out those
equations. [...] Where the really big mistakes are is in the
formulation of premises. [...] One of the
advantages of the mathematical medium -- or,
strictly speaking, of the mathematician's customary canons
of exposition of proof, whether in words or
in symbols -- is that we are forced to lay our cards
on the table so that all can see our premises.'' 

\smallskip

Returning to the large scale map of mathematical
provinces, Geometry/Algeb\-ra/Analysis, one should
find a place on it for (mathematical) Logic, with its modern
impersonation into the Theory of Algorithms and
Computer Science. There are compelling arguments
to consider it as a part of broadly conceived Algebra
({\it pace} Frege.)  And if one agrees on that,
Atiyah's insight about association of Algebra with Time,
becomes corroborated. In fact, the great shift
in the development of Logic in the 30s of the twentieth century 
occurred when Alan Turing used a physics
metaphor, ``Turing machine'', for the description
of an algorithmic computation. Before his work,
Logic was considered almost exclusively in
para--linguistic terms, as we did above.
Turing's vision of a finite automaton
moving in discrete steps along 
one--dimensional tape and writing/erasing bits
on it, and theorem about existence of a universal
machine of this type, stress exactly this
temporal aspect of all computations.
Even more important, the idea of computation as a physical process
not only helped create modern computers, but
also opened way to thinking in physical terms,
both in classical and quantum mode, about general laws of 
storing and processing information. 

\medskip

{\bf I.2. Objects of mathematical knowledge.} When we
study biology, we study living organisms. When we study astronomy, 
we study celestial bodies. When we study chemistry,
we study varieties of matter and of ways it can
transform itself.

\smallskip
 
We make observations and measurements of raw reality, 
we devise narrowly targeted experiments in a controlled environment
(not in astronomy however), and finally we
produce an explanatory paradigm, which becomes a current milestone
of science.

\smallskip

But what are we studying when we are doing mathematics?

\smallskip

A possible answer is this: {\it we are studying ideas
which can be handled as if they were real things.}
(P.~Davis and R.~Hersh call them ``mental objects with reproducible properties'').

\smallskip

Each such idea must be rigid enough in order to keep its
shape in any context it might be used. At the same time,
each such idea must have a rich potential of making connections
with other mathematical ideas. When an initial complex of
ideas is formed (historically, or pedagogically), connections
between them may acquire the status of mathematical objects
as well, thus forming the first level of a great hierarchy of abstractions.

\smallskip

At the very base of this hierarchy are mental images of things
themselves and ways of manipulating them. Miraculously,
it turns out that even very high level abstractions
can somehow reflect reality: knowledge of the world
discovered by physicists can be expressed only in the
language of mathematics.

\smallskip

Here are several basic examples.

\smallskip

{\bf I.2.1. Natural numbers.} This is arguably the oldest 
{\it proto--mathematical} idea. ``Rigidity'' of 1, 2, 3, ... is such that
first natural numbers acquire symbolic and religious meanings
in many cultures. Christian Trinity, or Buddhist Nirvana come to mind:
the latter evolved from Sanskrit {\it nir--dva--n--dva}
where {\it dva} means ``two'', and the whole expression
implies that the state of absolute blessedness is attained
through the extinction of individual existence
and becoming ``one'' with Universe. (These negative connotations of
the idea of ``two'' survive even in modern European languages
where it carries association with the idea
of ``doubt'' :
cf. Latin {\it dubius}, German {\it Zweifeln}, and 
Goethe's description of Mephistopheles). 

\smallskip

Natural number is also a {\it proto--physical} idea:
counting material objects (and later immaterial objects as well,
such as days and nights) is the first instance
of {\it measurement,} cf. below. 

\smallskip

Natural number becomes a {\it mathematical} idea
when: 

\smallskip

a) Ways of handling natural numbers as if they were 
things are devised: adding, multiplying. 

\smallskip

b) The first abstract features of the internal structure
of the totality of all natural numbers is discovered:
prime numbers, their infinity, existence and uniqueness
of prime decomposition.

\smallskip

These two discoveries were widely separated 
historically and geographically; arguably, 
culturally and philosophically as well.
Place--value system marks the origin of what we nowadays
call {\it applied mathematics}, primes mark the origin
of what used to be called {\it pure mathematics.}
Here are a few details.

\smallskip

At first, both  numbers and
ways of handling them are encoded by specific material objects:
fingers and other body parts, counting sticks, notches.
Notch is already a sign, not a proper thing, and it may
start signifying not 1, but 10 or 60, depending on where in the 
row of other symbols it is situated.
A way to the early great mathematical discovery,
that of place--value numeration system is open.
However, a consistent place--value system also requires a
sign for ``zero'', which came late and marked a new level
of mathematical abstraction.

\smallskip

An expressive summary in [B-BH] sketches the following
picture:

\smallskip 

``In c. 2074 BCE, king Shulgi organized a military reform 
in the Sumerian Empire, and the next year an administrative reform (seemingly introduced under the pretext of a state of emergency but soon made
permanent) enrolled the larger part of the working population
in quasi--servile labour crews and made overseer scribes
accountable for the performance of their crews, calculated in
abstract units worth $1/60$ of a working day (12 minutes)
and according to fixed norms. In the ensuing bookkeeping,
all work and output therefore had to be calculated
precisely and converted into these abstract units, which asked for
multiplications and divisions {\it en masse.} Therefore,
a place value system with base 60 was introduced for intermediate
calculations. Its functioning presupposed the use of tables of multiplication,
reciprocals and technical constants and the training for their
use in schools; the implementation of a system whose basic idea was
``in the air'' for some centuries therefore asked for decisions made at the level of the state and implemented with great force. Then as in many
later situations, only war provided the opportunity for such
social willpower.''

\smallskip

Primes, on the other hand, seem to spring off
from pure contemplation, as well as the idea
of a very concrete infinity, that of natural numbers
themselves, and that of prime numbers.

\smallskip

The proof of infinity of primes codified in Euclid's {\it Elements}
is a jewel of an early mathematical reasoning.
Let us recall it briefly in modern notation:
having a finite list of primes $p_1, \dots ,p_n$,
we can add one more prime to it by taking any
prime divisor of $p_1 \dots p_n +1$.

\smallskip

This is a perfect example of handling mathematical ideas
as if they were rigid material objects. And at this stage,
they are already
pure ideas bafflingly unrelated to any vestiges
of Sumerian or whatever {\it material notation}.
Looking at the modern decimal notation of a number,
one can easily tell whether it is even or divisible by 5, but not
whether it is prime. 
Generations of mathematicians after Euclid
marveled at an apparent randomness
with which primes pop up in the natural series.

\smallskip

Observation, controlled experimentation, and recently even
engineering of primes (producing and recognizing large
primes by computationally feasible algorithms, for security
applications) became a trademark of much of modern number
theory.

\medskip

{\bf 1.2.2. Real numbers and ``geometric algebra''.}
Integers resulted from counting, but other real numbers
came from geometry, as lengths and surfaces, volumes.
The discovery by Pythagoras of the incommensurability
of the diagonal of a square with its side
was at the same time the demonstration that
there were more ``magnitudes'' that ``numbers''.
Magnitudes were later to become real numbers. 

\smallskip
 
Arithmetical operations on integers evolved from
putting together sticks and notches to systematic
handling normalized notations in an ordered way.
Algebraic operations on reals evolved from drawing
and contemplating sketches which could intermittently
be plans of building sites or results of surveying, and renderings of
Euclidean circles, squares and angles.

\smallskip

Historians of mathematics in the twentieth century argued {\it pro}
and {\it contra} interpretation of a considerable part
of Greek mathematics as ``geometric algebra''.
One example of it is a sketch of a large square
subdivided into four parts by two lines
parallel to the orthogonal sides so that two
of the parts are again squares. This sketch can
be read as an expression and a proof
of the algebraic identity $(a+b)^2=a^2+b^2+2ab.$

\smallskip

Our modernizing perspective suggests a more general
consideration of several modes of mental processes,
in particular those related to mathematics.
The following two are the basic ones:

\smallskip

a) Conscious handling of a finite and discrete
symbolic system, with explicitly prescribed laws of formation
of meaningful strings of symbols, constructing new strings,
and less explicit rules of deciding which strings
are ``interesting'' (left brain, linguistic,
algebraic activity).

\smallskip

b) Largely subconscious handling of visual images,
with implicit reliance upon statistics of past experience,
estimating probabilities of future outcomes,
but also judging balance, harmony,
symmetry (right brain, visual arts and music, geometry).

\smallskip

Mental processes of mathematicians doing research
must combine these two
modes in many sophisticated ways. This is not an easy 
task, in particular because information processing rates
are so astonishingly different, of the order
10 bit/sec for conscious symbolic processing,
and $10^7$ bit/sec for subconscious visual processing
(cf. [N{\o}]).    

\smallskip

Probably because of inner tension created by this (and other)
discrepancies, they tend to be viewed emotionally,
as an embodiment of values -- cold intellect against 
warm feeling, bare logic vs. penetrating intuition.
See beautiful articles by David Mumford [Mu1] and [Mu2]
who eloquently defends statistics against logic,
but invokes mathematical statistics, which is built,
as any mathematical discipline, in an extremely logical way. 
 
\smallskip

Returning to real numbers and the ``geometric algebra'' of the Greeks,
we recognize in it a sample of right brain treatment
of a subject which later historically evolved 
into something dominated by the left brain. Or, as
Mumford puts it, modern algebra is a grammar
of actions with objects which are inherently geometric,
and Greek algebra is an early compendium of such actions.

\smallskip

Perhaps the continuity of Greek geometric thinking as a cognitive
phenomenon can be traced not only in modern geometry
but also in theoretical physics. The last decades have seen 
such a vigorous input of insights, 
conjectures, and sophisticated constructions, from
physics to mathematics, that an expression 
``physical mathematics'' was coined. The theoretical thinking
underlying the creative use of Feynman's path integral,
strikes us by the richness of results constructed on
a foundation which is mathematically shaky by
any standards. This can be considered as an additional
justification of the notion that ``geometric algebra''
was a reality, and not only our reconstruction of it.

\medskip

{\bf I.2.3. $e^{\pi i}= -1$: a tale of three numbers.}
Arguably, Euler's formula $e^{\pi i}= -1$ is the most
beautiful single formula in all mathematics.

\smallskip

It combines in a highly unexpected way three
(or four, if one counts $-1$ separately) constants
that were discovered in various epochs,  and 
emanate an aura of very different motivations.

\smallskip

Very briefly, $\pi = 3,1415926 ...$ is a legacy of the Greeks (again).
Even its existence as a real number, that is (the length of) 
a line segment, or surface of a square, is not something that
can be grasped without an additional mental effort.
The problem of ``squaring the circle'' is 
not just the next geometric problem, but a 
legitimacy test, with an uncertain outcome.

\smallskip

By contrast, $e = 2, 7128128 ...$ emerged  in the
already mature, if not fully developed Western mathematics
(mid--seventeenth century).
It is a combined theoretical by--product of the
invention of logarithm tables
as a tool of optimization of numerical algorithms 
(addition replacing multiplication) and the problem of
``squaring the hyperbole''. None of the
classical geometric constructions led to $e$
and none suggested any relation between $e$ and $\pi$.

\smallskip

Finally, the introduction of $i=\sqrt{-1}$, 
an ``imaginary'' number,
a monstrosity for many contemporaries,
was literally imposed on Cardano by the 
formulas for roots of a cubic equation
expressed in radicals. When all three roots
are real, formulas required complex numbers in
intermediate calculations.  

\smallskip

Euler's formula is a remarkable example of 
``infinite'' identities of which he (and later Srinivasa Ramanujan)
was a great practitioner. In fact,
$e^{\pi i}=-1$ is a particular case of the series $e^{ix}=
\sum_{n=0}^{\infty} (ix)^n/n!$ which gives
a more general expression $e^{ix}=\roman{cos}\,x+ i\,\roman{sin}\,x$.

\smallskip

Further progress in our understanding of real numbers and
theory of limits relegated the Euler and Ramanujan
great skills of dealing with ``infinite identities''
to backstage. G. Hardy, describing
Ramanujan's mathematical psyche, was at a loss
trying to interiorize it. This story does tell something
about the {\it logic vs statistics} dichotomy, but I cannot
pinpoint even a tentative statement.

\smallskip

As a totally unrelated development, 
$e^{ix}=\roman{cos}\,x+ i\,\roman{sin}\,x$ turned out
to be at the base of an adequate description 
of one of the most important and unexpected discoveries of the physics
of twentieth century: quantum probability amplitudes,
their wave--like behaviour, and quantum interference.

\medskip

{\bf I.2.4. Cantorian set: the ultimate mathematical object.}
In the original description by Cantor,

\smallskip

{\it Unter einer  `Menge'
verstehen wir jede $\underline{Zusammenfassung}$ $M$ von bestimmten
wohlunterschiedenen Objekten $m$ unserer Anschauung
oder unseres Denkens (welche die `Elemente' von $M$ genannt werden)
$\underline{zu\ einem\ Ganzen}$.}

\smallskip

``By a `set' we mean any collection $M$ into a whole of definite,
distinct objects $m$ (called the `elements' of $M$)
of our perception  or our thought.''

\smallskip

German syntax allows Cantor to mirror the meaning of the sentence
in its structure: {\it Objekten $m$ unserer Anschauung etc} are packed  
between the opening bracket {\it {Zusammenfassung}} 
and the closing bracket {\it zu einem Ganzen}. 

\smallskip

Contemplating this definition for the first time,
it is difficult to imagine what kind of mathematics 
or, for that matter, what kind of mental activity at all,
can be performed with such meager means.
In fact, it is precisely this parsimony which allowed Cantor
to invent his ``diagonal process'', 
to compare infinities as if they were physical objects,
and to discover that the infinity of real numbers is
strictly larger than that of integers.

\smallskip

Simultaneously, Cantor's intuition underlies most
of foundational work in the mathematics of the
twentieth century: it is either vigorously refuted
by logicists of various vintages, or works
as a great unification project, in both guises of Set Theory
and its successor, Category Theory.

\medskip

{\bf I.2.5. ``All men are mortal, Kai is a man ...'': from syllogisms
to software.} Aristotle codified elementary forms of statements and
basic rules of logical deductions. The analogies
between them and elementary arithmetics
were perceived early,  but made precise late;
we recognize Boole's role in this development.
Philosophers of science disagreed about hierarchical
relationships between the two. Frege, for example,
insisted that arithmetic was a part of logic.

\smallskip

The 20th century has seen a sophisticated fusion
of both realms when in the thirties G\"odel, Tarski and Church
produced mathematical models of mathematical reasoning
going far beyond the combinatorics of finite texts.
One of the important tools was the idea, going back to Leibniz,
that one can use a computable enumeration of all texts by
integers allowing to replace logical deductions by
arithmetical operations.

\smallskip

Tarski modeled truth as ``truth in all interpretations'',
and found out that the set of (numbers of) arithmetical
truths  cannot be expressed by an arithmetical formula.
Infinitarity of Tarski's notion of truth is connected 
with the fact that logical formulas are allowed to contain
quantifiers ``for all'' and ``there exists'',
so interpretation of a finite formula involves
potentially infinite sequence of verifications.  

\smallskip

G\"odel, using a similar trick, demonstrated that the set
of arithmetical truths deducible from any finite 
system of axioms and deduction rules cannot coincide
with the set of all true formulas. Self--referentiality
was an essential common feature of both proofs.

\smallskip

Among other things, G\"odel and Tarski  showed that 
the basic hierarchical relation
is that between a language and a metalanguage.
Moreover, only their interrelation and not absolute
status is objective. One can use logic to describe
arithmetics, and one can use arithmetics to discuss
logic. A skillful mixture of both levels
unambiguously shows inherent restrictions
of pure logic as a cognitive tool, even when
it is applied ``only'' to pure logic itself.

\smallskip

Turing and Church during the same decade
analyzed the idea of ``computability'',
which had a more arithmetic flavor from the start.
Alan Turing made a decisive step by substituting a physical
image (Turing machine) in place of the traditional
linguistic embodiments for logic and computation
dominating both Tarski's and G\"odel's discourses. 
This was a great mental step
preparing the subsequent technological evolution:
the emergence of programmable electronic calculators.

\smallskip

Theoretically, both Church and Turing discovered that there existed
a ``final'' notion of computability embodied in
the universal recursive
function, or universal Turing machine. This was 
not a mathematical theorem, but rather a 
``physical discovery in a metaphysical realm'',
justified not by a proof but by the fact that 
all subsequent attempts to conceive an alternative version
led to an equivalent notion. A ``hidden'' (at least in popular
accounts) part of this discovery was the realization that
the correct definition of computability includes
elements of un--computability that cannot be avoided 
at any cost: a recursive function is generally not
everywhere defined, and we cannot decide
at which points it is defined and at which not.

\smallskip

Computers which are functioning now, embody
a technologically alienated form of these great insights. 

\medskip

{\bf I.3. Definitions/Theorems/Proofs.} I will briefly describe now
tangible traces of ``pure'' mathematics as a collective activity
of the contemporary professional community. I will stress not so much
organizational forms of this activity as external
reflection of the inner structure of the world of mathematical ideas.

\smallskip
 
Look at any contemporary paper in one of the leading
research journals like {\it Annals of Mathematics}
or {\it Inventiones mathematicae.} Typically,
it is subdivided into
reasonably short patches called Definitions,
Theorems (with Lemmas and Propositions as subspecies),
and Proofs, that can be considerably longer.
These are the basic structure blocks of a modern
mathematical exposition; frills like
motivation, examples and counterexamples,
discussion of special cases, etc., make it livelier.

\smallskip

This tradition of organizing mathematical knowledge
is inherited from the Greeks, especially Euclid's Elements.
The goal of a definition is to introduce a  mathematical
object. The goal of a theorem is to state some of its
properties, or interrelations between various objects.
The goal of a proof is to make such a statement
convincing by presenting a reasoning subdivided into small
steps each of which is justified as an ``elementary'' convincing
argument.

\smallskip

To put it simply, we first explain, what we are talking about,
and then explain why what we are saying is true
({\it pace} Bertrand Russell).

\smallskip

{\it Definitions.} The first point is epistemologically subtle and controversial,
because what we are talking about are extremely
specific mental images not present normally in an
untrained mind (what is a a real number? a
random variable? a group?). Presenting some basic objects
above, I used narrative devices to make them look
more graphic or vivid, but gave no real definitions
in the technical sense of the word.  

\smallskip

Euclid's definitions 
usually consist of a mixture of explanations
involving visual images, and ``axioms'' involving
some idealized properties that we want to impose on them.

\smallskip

In contemporary mathematics, one can more or less
explicitly restrict oneself to the basic
mental image of a Cantorian ``set'', and 
a limited inventory of properties of
sets and constructions of new sets from given ones.
Each of our Definitions then can be conceived as a
standardized description of a certain structure,
consisting of sets, their subsets etc.
This is a viewpoint that was developed by the Bourbaki 
group and which proved to be an extremely influential,
convenient and widely accepted way of organizing mathematical 
knowledge. Inevitably, a backlash ensued,
aimed mostly at the value system supporting this
neo--Euclidean tradition, but its pragmatic
merits are indisputable. At the very least,
it enabled a much more efficient communication
between mathematicians coming from different fields.

\smallskip

If one adopts a form of Set Theory as a basis
for further constructions, only set--theoretic axioms
remain ``axioms'' in Euclid's sense, something like
intuitively obvious properties accepted without
further discussion (but see below), whereas the
axioms of real numbers or of plane geometry
become provable properties of explicitly constructed
set--theoretic objects.

\smallskip

Bourbaki in their multivolume treatment
of contemporary mathematics developed
this picture and added to it
a beautiful notion of ``structures--m\`eres''
(the issue [Sci] is dedicated to the history of the Bourbaki group).

\smallskip

In a broader framework, one can  argue
that mathematicians have developed a specific discursive behavior
which might be called ``culture of definitions''.
In this culture, many efforts are invested
into clarification of the content (semantics) of basic abstract
{\it notions} and syntax of their interrelationships, whereas the choice
of {\it words} (or even to a larger degree, {\it notations}) 
for these notions is a secondary matter and largely arbitrary convention,
dictated by convenience, aesthetic considerations, 
by desire to invoke appropriate connotations.
This can be compared with some habits of humanistic
discourse where such terms as {\it Dasein} or {\it diff\'erance}
are rigidly used as markers of a certain tradition,
without much fuss about their meaning.

\medskip

{\bf I.4. Problems/Conjectures/Research Programs.}
From time to time, a paper appears which solves,
or at least presents in a new light,
some great problem, or conjecture, which was with us for
the last decades, or even centuries, and resisted
many efforts. Fermat's Last Theorem (proved by Andrew Wiles),
the Poincar\'e Conjecture, the Riemann Hypothesis,
the P/NP--problem these days even  make newspapers headlines.

\smallskip

David Hilbert composed his talk at the second 
(millennium) International Congress of Mathematicians
in Paris on August 8, 1900, as a discussion of ten outstanding
mathematical problems which formed a part of his list
of 23 problems compiled in the published version. 
One can argue about their
comparative merit in pure scientific terms,
but certainly they played a considerable role in focusing
efforts of mathematicians on well defined directions,
and providing clear tasks and motivation for young
researchers. 

\smallskip

Whereas a problem (a {\it yes/no} question) is basically a guess
about validity or otherwise of a certain statement
(like Goldbach's problem: every even number
$\ge 4$ is a sum of two primes), a Research Program is 
an outline of a broad vision, a map of a landscape
some regions of which are thoroughly investigated, 
whereas other parts are guessed on the base of analogies,
experimentation with simple special cases, etc.

\smallskip

The distinction between the two is not absolute.
Problem Number one, the Continuum Hypothesis, which
in the epoch of Cantor and Hilbert looked like a
{\it yes/no} question, generated a vast research program
which established, in particular, that neither of the
two answers is deducible within the generally accepted
axiomatic Set Theory. 

\smallskip

On the other hand, the explicit formulation 
of a research program can be a risky venture.
Problem Number 6 envisioned the axiomatization of physics.
In the next three decades or so physics 
completely changed its face.

\smallskip

Some of the most influential Research Programs of
the last decades were expressions of insights
into the complex structure of Platonian reality.
A.~Weil guessed the existence of cohomology theories
for algebraic manifolds in finite characteristics.
Grothendieck constructed them, thus forever changing our
understanding of the relationships between continuous and discrete.

\smallskip

When Poincar\'e said that there are no solved problems,
there are only problems which are more or less solved,
he was implying that any question
formulated in a {\it yes/no} fashion is
an expression of narrow--mindedness. 

\smallskip

The dawning of the twenty first century was marked by the publication 
by the Clay Institute of the list of
Millenium Problems. There are exactly seven of them,
and they are all {\it yes/no} questions.  For the
first time a computer science--generated problem appears:
the famous P/NP conjecture. 
Besides, Clay Problems come with a price tag: USD $10^6$
for a solution of any one of them. Obviously,
free market forces played no role in this pricing policy.

\bigskip

\centerline{\bf II. Mathematics as a Cognitive Tool}

\medskip

{\bf II.1. Some history.} Old texts that are considered as sources
for history of mathematics show that it
started as a specific activity answering the needs of commerce and
of state, servicing large communal works and warfare:
cf the excerpt above about Sumero--Babylonian
administrative reform. 

\smallskip

As another example, turn to the Chinese book
``The nine chapters on mathematical procedures''
compiled during the Han dynasty around the beginning of our era.
We rely here upon the report of K.~Chemla at the Berlin ICM 1998,
[Che]. The book generally is a sequence of problems
and of their solutions which can be read
as special cases of rather general algorithms so that a 
structurally similar
problem with other values of parameters could be solved as well.
According to Chemla, problems ``regularly invoke
concrete questions with which the bureaucracy of the Han dynasty was faced, and, more precisely, questions that
were the responsibility of the ``Grand Minister of Agriculture'' 
({\it dasinong}), such as renumerating civil servants, managing granaries or enacting standard grain measures. Moreover, the sixth of 
{\it The nine chapters} takes its name from an economic
measure actually advocated by a Grand Minister of Agriculture, Sang Hongyang (152--82 B.C.E), to levy taxes in a fair way, a program for which 
the Classic provides mathematical procedures.''

\smallskip

Yet another description of the preoccupations
of Chinese mathematicians is given in [Qu]:

\smallskip

``In the long history of the Chinese empire, mathematical
astronomy was the only subject of the exact sciences that attracted
great attention from rulers. In every dynasty, the 
royal observatory was an indispensable part of the state. Three
kinds of expert -- mathematicians, astronomers and astrologers --
were employed as professional scientists by the emperor.
Those who were called mathematicians took charge of 
establishing the algorithms of the calendar--making systems. 
Most mathematicians
were trained as calendar--makers. [...]

Calendar--makers were required to maintain a high degree 
of precision in prediction. Ceaseless efforts to improve 
numerical methods were made in order to
guarantee the precision required for astronomical observation [7].
It was neither necessary nor possible that a geometric model
could replace the numerical method, which occupied the principal
position in Chinese calendar--making system. [...] As a subject
closely related to numerical method, algebra, rather than
geometry, became the most developed field of mathematics
in ancient China.''

\smallskip

Western tradition goes back to Greece. 
According to Turnbull [Tu], we owe the word 
``mathematics'' and the subdivision of mathematics
into Arithmetic and Geometry to Pythagoras (569 -- 500 BC).
More precisely, Arithmetic (and  Music)
studies the discrete, whereas Geometry and Astronomy
study the continued. The secondary dichotomy Geometry/Astronomy
reflects the dichotomy The stable/The moving.

\smallskip

With small modifications, this classification was at the origin 
of the medieval ``Quadrivium of knowledge'', and 
Michael Atiyah's overall view of mathematics still bears
distinctive traces of it.

\smallskip

Plato (429--348 BC) in {\it Republic}, Book VII, 525c, explains why
the study of arithmetic is essential
for an enlightened statesman:

\smallskip

``Then this is a kind of knowledge, Glaucon,
which legislation may fitly prescribe; and we must endeavour to persuade those who are prescribed to be the principal men of our State to go and learn arithmetic, not as amateurs, but they must carry on the study until they see the nature of numbers with the mind only; nor again, like merchants or retail-traders, with a view to buying or selling, but for the sake of their military use, and of the soul herself; and because this will be the easiest way for her to pass from becoming to truth and being.''

\smallskip

With gradual emergence of ``pure mathematics'', return to
practical needs began to be classified as applications.
The opposition pure/applied mathematics as we know it now
certainly has already crystallized by the beginning of the nineteenth
century. 
In France, Gergonne was publishing the {\it Annales de math\'ematiques pures et appliqu\'ees} which ran from 1810 to 1833. In Germany
Crelle founded in 1826 the {\it Journal f\"ur die reine und angewandte
Mathematik.}

\medskip

{\bf II.2. Cognitive tools of mathematics.} In order to understand
{\it how} mathematics is applied to the understanding of real world,
it is convenient to subdivide it into the following three
modes of functioning: {\it model, theory, metaphor.}

\smallskip

A mathematical {\it model} describes a certain range of phenomena qualitatively
or quantitatively but feels uneasy pretending to be something more.

\smallskip

From Ptolemy's epicycles (describing planetary motions, ca 150) 
to the Standard Model (describing interactions of elementary particles, 
ca 1960), quantitative models cling to the observable reality
by adjusting numerical values of sometimes dozens of 
free parameters ($\ge 20$ for the Standard Model). Such models 
can be remarkably precise. 

\smallskip

Qualitative models offer insights into {\it stability/instability},
{\it attractors} which are limiting states tending to occur independently
of initial conditions, {\it critical phenomena} in complex
systems which happen when the system crosses a boundary
between two phase states, or two basins
of different attractors. A recent report [KGSIPW] is dedicated
to predicting of surge of homicides
in Los Angeles, using as  methodology the pattern recognition of
infrequent events. Result: ``We have found that
the upward turn of the homicide rate is preceded within 11 months
by a specific pattern of the crime statistics:
{\it both burglaries and assaults simultaneously escalate, while robberies
and homicides decline.} Both changes, the escalation and
and the decline, are not monotonic,
but rather occur sporadically, each lasting some
2--6 months.''

\smallskip

The age of computers has seen a proliferation of models,
which are now produced on an industrial scale
and solved numerically. A perceptive essay by 
R.~M.~Solow ([Sol], written in 1997) argues that modern
mainstream economics is mainly concerned with model--building.

\smallskip

Models are often used as ``black boxes''
with hidden computerized input
procedures, and oracular outputs prescribing behavior of 
human users, e.~g. in financial transactions.

\smallskip

What distinguishes a (mathematically formulated physical) {\it theory} 
from a model is primarily its 
higher aspirations. A modern physical
theory generally purports that it would describe the world
with absolute precision if only it (the world)
consisted of some restricted variety of
stuff: massive point particles obeying only the
law of gravity; electromagnetic field in a vacuum;
and the like. In Newton's law for the force $\frac{Gm}{r^2}$
acting on a point in the central gravity field,
$Gm$ and $r$ might be concessions to measurable reality, but 
$2$ in $r^2$ is a rock solid theoretical $2$, not some $2,000000003...$,
whatever experimentalists might measure to the contrary.
A good quantitative theory can be very useful
in engineering: a machine is an artificial fragment
of the universe where only a few physical laws are allowed to dominate
in a well isolated material environment. In this function,
the theory supplies a model. 

\smallskip

A recurrent driving force generating theories is a concept
of a reality beyond and above the material world,
reality which may be grasped only by mathematical tools.
From Plato's solids to Galileo's ``language of nature''
to quantum superstrings, this psychological attitude
can be traced sometimes even if it conflicts with the
explicit philosophical
positions of the researchers.    

\smallskip

A (mathematical) {\it metaphor}, when it aspires
to be a cognitive tool, postulates that some 
complex range of phenomena might be compared to
a mathematical construction. The most recent
mathematical metaphor I have in mind is
Artificial Intelligence (AI). On the one hand, AI
is a body of knowledge related to computers and
a new, technologically created reality, consisting of
hardware, software, Internet etc. On the other hand,
it is a potential model of functioning of biological
brains and minds. In its entirety, it has not reached
the status of a model:
we have no systematic, coherent and extensive
list of correspondences between chips and neurons,
computer algorithms and brain algorithms.
But we can and do use our extensive knowledge of 
algorithms and computers
(because they were created by us) to generate
educated guesses about structure and function of the central
neural system: see [Mu1] and [Mu2].

\smallskip

A mathematical theory is an invitation 
to build applicable models. A mathematical metaphor
is an invitation to ponder upon what we know.
Susan Sontag's essay about (mis)uses of the ``illness''
metaphor in [So] is a useful warning.

\smallskip

Of course, the subdivision I have just sketched 
is not rigid or absolute. 
Statistical studies in social sciences often
vacillate between models and metaphors.
With a paradigm change,
scientific theories are relegated to the status
of outdated models. But for the sake of our exposition,
it is a convenient way to organize synchronic and historical data.

\smallskip

I will now give some more details about these cognitive tools,
stressing models and related structures.

\medskip

{\bf II.3. Models.} One can analyze the creation and functioning
of a mathematical model by contemplating the following stages
inherent in any systematic study of quantifiable observations.

\smallskip

i) Choose a list of observables.

\smallskip

ii) Devise a method of measurement: assigning numerical
values to observables. Often this is preceded by a more or less explicit
ordering of these values along an axis (``more -- less'' relation);
then measurement is expected to be consistent with ordering.

\smallskip

iii) Guess the law(s) governing the distribution of
observables in the resulting, generally multidimensional,
configuration space.  The laws can be probabilistic or exact.
Equilibrium states can be especially interesting; they are often
characterized as stationary points of an appropriate functional
defined on the whole configuration space.
If time is involved, differential equations for evolution
enter the game.  

\medskip

Regarding the idea of ``axis'', one should mention
its interesting and general cultural connotations
expounded by Karl Jaspers. Jaspers postulated
a transition period to modernity around 500 BC,
an ``axial time'' when a new human mentality emerged 
based on the opposition between immanence and transcendence.
For us relevant here is the image of oppositions
as opposite orientations of one and the same axis,
and the idea of freedom as a freedom of choice between
two incompatible alternatives. This is also the imagery
behind the standard physical expression ``degrees of freedom'',
which is now almost lost, as usually happens to images
when they become terms.

\smallskip

The idea of measurement, which is the base
of modern science, is so crucial that it is sometimes
uncritically accepted in model--building.
It is important to keep in mind its
restrictions. 

\smallskip

In the quantum mode of description of the microworld,
a ``measurement'' is a very specific interaction 
which produces a random change
of the system state, rather than furnishing information
about this state.

\smallskip

In economics, money serves as the universal axis upon which
``prices'' of whatever are situated.
``Measurement'' is purportedly a function of market forces.

\smallskip

The core intrinsic contradiction of the market
metaphor (including the outrageous ``free market of ideas'') 
is this: we are projecting the multidimensional
world of incomparable and incompatible degrees of freedom to the
one--dimensional world of money prices. As a matter of principle,
one cannot make 
it compatible with even basic order relations on these axes,
much less compatible with non--existent or 
incomparable values of different kinds.

\smallskip

In this respect, the most oxymoronic use of the market metaphor
is furnished by the expression ``free market of ideas''.

\smallskip

Only one idea is on sale
at this market: that of ``free market''.

\smallskip

{\bf II.3.1. A brief glossary of measurement.} 
A general remark about measurement: for each ``axis''
we will be considering, the history of measurements starts with
the stage of ``human scale'' and involves direct manipulation
with material objects. Gradually it evolves to much larger and much smaller
scales, and in order to deal with the new challenges
posed by this evolution, more and more mathematics 
is created and used.

\smallskip

COUNTING. We suggest to the reader to reread the subsection
on Natural Numbers above as a glimpse into the history
of counting (and accounting). It shows clearly how the transition 
from counting small quantities of objects (``human scale'')
to the scale of state economy stimulated the creation
and codification of a place--value notation. 

\smallskip

Skipping other interesting developments,
we must briefly mention what Georg Cantor
justifiably considered as his finest achievement:
counting ``infinities'' and the discovery
that there is an infinite scale of infinities of 
growing orders of magnitude.

\smallskip

His central argument is structurally very similar to the
Euclid's proof that there are infinitely many primes:
if we have a finite or infinite set $X$, then the set of all its subsets
$P(X)$ has a strictly larger cardinality. This is established
by Cantor's famous ``diagonal'' reasoning.

\smallskip

Cantor's theory of infinite sets produces an incredible
extension of both aspects of natural numbers: each number measures 
``a quantity'', and they are ordered by the relation
``$x$ is larger than $y$''. Infinities, respectively,
are ``cardinals'' (measure of infinity) and ``ordinals''
which are points on the ordered axis of growing infinities.

\smallskip

The mysteries of Cantor's scale led to a series of unsolved
(and to a considerable degree unsolvable) problems, and became
the central point of
many epistemological and foundational discussions in the twentieth century.
The controversies and bitter arguments about the legitimacy
of his mental constructions made
the crowning achievement of his life
also the source of a sequence of nervous breakdowns and depressions
which finally killed him as the world war I was slowly grinding
the last remnants of Enlightenment's belief in reason.   

\smallskip 

SPACE AND TIME. Human scale measurements of length 
must have been inextricably related to those of plots, 
and motivated by agriculture and building.
A stick with two notches, or a piece of string, could be used
in order to transport a measure of length from one place to another.
 
\smallskip

Euclid's basic abstraction: an infinitely rigid
and infinitely divisible plane, with its hidden symmetry
group of translations and rotations, with its
points having no size, lines stretching uninterrupted in
two directions, perfect circles and triangles,
must have been a refined mental image of the ancient
geodesy. Euclid's space geometry arguably was even closer to the
observable world, and it is remarkable, that he
systematically produced and studied abstractions of
two--, one--, and zero-- dimensional objects as well.

\smallskip

Pythagoras's theorem was beautifully related to arithmetic
in the practice of Egyptian builders: the formula $3^2+4^2=5^2$
could be transported into a  prescription for producing
a right angle with the help of a string with uniformly
distanced knots on it.

\smallskip

When  Eratosphene of Alexandria (ca 200 BC) devised his
method for producing the first really large scale
scientific length measurement, that of the size of the Earth,
he used the whole potential of Euclid's geometry
with great skill. He observed that at noon on the day of 
summer solstice at Syene the sun was exactly at the zenith
since it shined down a deep well.
And at the same time at Alexandria the distance
of the sun to the zenith was one fiftieth of the circumference.
Two additional pieces of observational data were
used. First, that the distance between Syene and Alexandria,
which was taken to be 5000 Greek stades
(this is also a large scale measurement, probably, based upon the
time needed to cover this distance). Second, the assumption that
Syene and Alexandria  lie on the same meridian.

\smallskip

The remaining part of Eratosphene's measurement method 
is based upon a theoretical model.
Earth is supposed to be round, and Sun to be at an essentially
infinite distance from its center, so that the lines of sight
from Syene and Alexandria to the Sun are parallel.

\smallskip

Then an easy Euclidean argument applied to the cross--section
of the Earth and outer space passing through Syene, Alexandria, and the Sun,
shows that the 
distance between Syene and Alexandria must be one fiftieth
of the Earth circumference, which gives for the latter
the value 250000 stades. (According to modern evaluation
of Greek stade, this is a pretty good approximation.)

\smallskip

Implicit in this argument is an extended symmetry group
of the Euclidean plane including, with translations and rotations, also
rescalings: changing
all lengths simultaneously in the same proportion.
The practical embodiment of this idea, that of
{\it a map}, was crucial for an immense amount of human activities,
including geographical discoveries all over the globe. 

 \smallskip

The attentive reader has remarked already that
time measurements crept into this description
(based upon a book of Cleomedes ``De motu circulari corporum caelestium'',
middle of the first century BC). In fact, how do we know
that we are looking at the position of sun {\it at the same moment}
in Alexandria and Syene, distanced by 5000 stades? 

\smallskip

The earliest human scale time measurement
were connected with periodical cycles of day/night
and approximate position of sun on the sky. Sky dials,
referred to by Cleomedes and Eratosthenes, 
translate time measurements into space measurements.

\smallskip 

The next large scale measurements of time are
related to  the seasons of the year and periodicity of religious
events required in the community. 
Here to achieve the necessary precision, mathematical observational 
astronomy is needed. It is used first to register
irregularities in the periodicity
of year, so basically in the movement of Earth in the solar system.
The mathematics which is used here involves numerical calculation
based on interpolation methods. 

\smallskip

Next level of large scale: chronology of ``historical time''. 
This proved to be a rather un--mathematical endeavor. 

\smallskip

Geological and evolutionary time returns us to science:
the evolution of Earth structures and of life is traced
on the background of a well developed understanding of physical time 
which is highly mathematicized; however, the 
changes are so gradual and the evidence so scattered that 
precision of measurements ceases to be
accessible or essential. Besides the plethora of
observational data, brilliant guesses, and very elementary accompanying reasoning,
one small piece of mathematics becomes essential for dating:
the idea that radioactive decay leaves remnants of the decaying substance
whose quantity diminishes exponentially with time.
One very original version of this idea was used in ``glottochronology'':
the dating of proto--states of living languages which were reconstructed
using methods of comparative linguistics.

\smallskip

The sheer span of geological and evolutionary time when it was
first recognized and scientifically elaborated presented
a great challenge to the dogmata of (Christian) faith:
discrepancy with the postulated age of the World since the time of 
Creation became gaping.

\smallskip

Time measurements at a {\it small scale} become
possible with invention of clocks. Sundials
use relative regularity of visible solar motion
and subdivide daytime into smaller parts. Water and sand clocks
measure fixed stretches of time. This uses  the idea of reproducibility
of some well controlled physical processes. Mechanical clocks
add to this artificial creation of periodic
processes. Modern atomic clocks use subtle
enhancing methods for exploiting natural periodic
processes on a microscale.

\smallskip

Still, time remains a mystery, because we 
cannot freely move in it as we do in space, 
we are dragged to who knows where, and St Augustine reminds us
about this perennial, un--scientific torment:
``I know that I am measuring time. But I am not
measuring the future, for it is not yet; and I am not
measuring the present because it is extended
by no length; and I am not measuring the past
because  it no longer is. What is it, therefore,
that I am measuring? '' ({\it Confessions}, Book XI, XXVI.33).

\smallskip

CHANCE, PROBABILITY, FINANCE. Connotations
of the words ``chance'' and ``probability'' 
in the ordinary speech do not have much in common
with mathematical probability: see [Cha] for an 
interesting analysis of
semantics of related words in several ancient
and modern European languages.
Basically, they invoke the idea of human confidence 
(or otherwise) in an uncertain situation.

\smallskip

Measurements of probability, and mathematical handling
of the results,  refer {\it not} to the confidence itself
which is a psychological factor, but to
objective numerical characteristics of reality,
initially closely related to count.

\smallskip

If a pack contains 52 cards and they are well
shuffled, the probability to pick the queen of spades
is 1/52. Elementary but interesting mathematics  
enters when one starts calculating 
probabilities of various combinations (``good hands'').
Implicitly, such calculations involve
the idea of symmetry group: we not only count
the number of cards in the pack,
or number of good hands among all possible, but
assume that each one is equally probable
if the game is fair. 

\smallskip

The mathematics of gambling
was one source of probability theory, 
while another was the statistics of banking,
commerce, taxation etc. Frequencies of various
occurences and their stability led to
the notion of empirical probability
and to the more or less explicit idea of
``hidden gambling'', the unobservable realm
of causes which produced observable frequences
with sufficient regularity in order to fit into
a mathematical theory. The modern definition of a
probability space is an axiomatization of such an image.

\smallskip

Money started as a measure of value and made a 
crucial transition to the world of probability with the
crystallization of credit as a main function 
of a bank system. 

\smallskip

The etymology of the word ``credit''
again refers to  the idea of human confidence.
The emergent ``culture of finance'', according to
the astute analysis of Mary Poovey in [Po2],
drastically differs from an economy of production
``which generates profit by turning labor power into products
that are priced and and exchanged in the market''.
Finance generates profit, in particular,
``through placing complex wagers that future prices will raise or fall''
([Po2], p. 27), that is, through pure gambling.
The scale of this gambling is staggering,
and the incredible mixture of real and virtual
worlds in the culture of finance is explosive. 

\medskip 

INFORMATION AND COMPLEXITY. This is an example of a quite sophisticated
and contemporary measurement paradigm.

\smallskip

As with ``chance'' and ``probability'', the term 
{\it quantity} of information,
which became one of the important theoretical notions in the
second half of the twentieth century after the works of Claude Shannon
and Andrei Kolmogorov, has somewhat misleading
connotations. Roughly speaking, the quantity of
information is measured simply by the length
of a text needed to convey it. 

\smallskip

In the everyday usage, this measure seems to be rather 
irrelevant, first, and disorienting, second.
We need to know whether information is {\it important} and
{\it reliable}: these are qualitative rather than quantitative
characteristics. Moreover, importance
is a function of cultural, scientific, or political context. 
And in any case, it seems preposterous to
measure the information content of ``War and Peace'' by its 
sheer volume. 

\smallskip

However, quantity of information becomes central if we are
handling information without bothering about
its content or reliability (but paying attention to
security), which is the business of the media and
communication industry. The total size of texts
transmitted daily by Internet, mass media and phone
services is astounding and far beyond the limits of what
we called ``human scale''.

\smallskip

Shannon's basic ideas about measuring quantity of information
can be briefly explained as follows. Imagine
first that the information you want to transmit
is simply the answer ``yes'' or ``no'' to a question of
your correspondent. For this, it is not even necessary
to use words of any natural language: simply transmit
1 for ``yes'' and 0 for ``no''. This is one bit of
information. Suppose now that you want to transmit
a more complex data and need a text containing $N$ bits.
Then the quantity of information you transmit
is at least bounded from above by $N$, but
how do you know that you cannot use a shorter text
to do the same job? In fact, there exist {\it systematic} methods
of compressing the raw data, and they were made explicit
by Shannon. The most universal of them starts with the assumption
that in the pool of texts you might be wanting to transmit
not all are equally probable. In this case
you might change encoding in such a way that
the more probable texts will get shorter codes
than less probable ones, and thus save on the volume
of transmission, at least in average.

\smallskip

Here is how one can do it in order to encode
texts in a natural language. Since there are about 
$30$ letters of alphabet, and $2^5=32$, one needs 
5 bits to encode each one,
and thus to get a text whose bit--length is about 5 times
its letter--length. But some letters statistically 
are used much more often than others, so one can try to
encode them by shorter bit sequences. This leads to
an optimization problem that can be explicitly solved,
and the resulting length of an average compressed text
can be calculated. This is essentially the definition
of Shannon's and Kolmogorov's {\it entropy}.

\smallskip

Using the statistical paradigm of  measurement, the creators of 
Google found an imaginative solution
for the problem of assigning numerical measure
to the {\it relevance} of information as well.
Roughly speaking, a search request makes Google  produce
a list of pages containing a given word or expression.
Typically, the number of such pages is very large, 
and they must be presented
in the order of decreasing importance/relevance.
How does Google calculate this order?

\smallskip

Each page has hypertext links to other pages. One can model the whole
set of pages on the Web by the vertices of an oriented graph whose edges are
links. One can assume in the first approximation that importance of a
page can be measured by the number of links pointing to it. But this
proposal can be improved upon, by noting that all links are not equal: a
link from an important page has proportionately more weight, and a link
from a page that links to many other pages has proportionately less
weight. This leads to an ostensibly circular definition (we omit a
couple of minor details): each page imparts its importance to the pages
it links to, divided equally between them; each page's importance is
what it receives from all pages that link to it. However a classical
theorem due to A. Markov shows that this prescription is well defined.
It remains to calculate the values of importance and to range pages in
their decreasing order.

\smallskip

Let us now return to the
Shannon's optimal encoding/decoding procedures.
The reader has noticed that economy on transmission 
has its cost: encoding at the source and decoding
at the target of information.

\smallskip

What happens if we allow more complex encoding/decoding
procedures in order to achieve further degree of compression?

\smallskip

The following metaphor here might be helpful:
an encoded text at the source is essentially {\it a program} $P$
for obtaining the decoded text $Q$ at the target.
Let us now allow to transmit {\it arbitrary}
programs that will generate $Q$; perhaps we will be able
to choose the shortest one and to save resources.

\smallskip

A remarkable result due to Kolmogorov is that this is a well defined
notion: such shortest programs $P$ exist and their
length ({\it the Kolmogorov complexity of $Q$}) 
does not depend essentially  on the programming
method. In other words, there exists
{\it a totally objective measure of the quantity of
information contained in a given text $Q$.}

\smallskip

Bad news, however, crops up here: a) one cannot
systematically  reconstruct $P$ knowing $Q$ (unlike the case of Shannon entropy);
b) it may take a very long time decoding $Q$
from $P$ even if $P$ is known and short. A very simple example:
if $Q$ is a sequence of exactly $10^{10^{10}}$\ 1's,
one can transmit this sentence, and let the addressee
bother with the boring task of printing $10^{10^{10}}$\ 1's out.

\smallskip

This means that Kolmogorov's complexity, a piece of beautiful
and highly sophisticated (although ``elementary'') mathematics,
 is not a practical 
measure of quantity of information. However,
it can be used as a powerful metaphor elucidating
various strengths and weaknesses of the modern information society.

\smallskip

It allows us to recognize one essential way
in which scientific (but also everyday life)
information used to be encoded. The basic physical ``laws of nature''
(Newton's $F=ma$, Einstein's $E=mc^2$, the Schr\"odinger equation etc. )
are very compressed programs for obtaining relevant
information in concrete situations. Their Kolmogorov complexity
is clearly of human size, they bear names of humans associated with their
discovery, and their full information content is totally 
accessible to a single mind of a researcher or a student.

\smallskip

Nowadays, such endeavors as
the Human Genome projects provide us with huge
quantities of scientific data whose volume in any compressed 
form highly exceeds the capability of any single mind.
Arguably, similar databases that will be created for 
understanding the central nervous system (brain)
will present the same challenge, having Kolmogorov
complexity of comparable size with their volume. 

\smallskip

Thus, we are already studying those domains of material world
whose descriptions
have much higher information content (Kolmogorov complexity)
than the ones that constituted the object of classical science. 
Without computers,
neither the collective memory of observational data 
nor their processing would be feasible.

\smallskip

What will happen when the total 
essential new scientific ``knowledge'' {\it and}
its handling will have to be relegated to large computer 
databases and nets?

\bigskip

\centerline{\bf III.  Mathematical Sciences and Human Values}

\medskip

{\bf III.1. Introduction.}
Commenting on the fragments of the Rhind papyrus, a handbook
of Egyptian mathematics written about 1700 BC, 
the editor of the whole anthology [WM] James R.~Newman
writes (vol. I, p. 178, published in 1956) :

\smallskip

 ``It seems to me that a sound appraisal of
Egyptian mathematics depends upon a much broader and deeper 
understanding
of human culture than either Egyptologists or historians
of science are wont to recognize. As to the question
how Egyptian mathematics compares with Babylonian
or Mesopotamian or Greek mathematics,
the answer is comparatively easy and comparatively unimportant. 
What is more to the point is to understand why the Egyptians produced 
their particular kind of mathematics, to what extent it
offers a culture clue, how it can be related
to their social and political institutions, to their religious
beliefs, their economic practices, their habits of daily
living. It is only in these terms that their
mathematics can be judged fairly.''

\smallskip

By 1990, this became a widely accepted paradigm, and 
D'Ambrosio coined the term ``Ethnomathematics'' for it
(cf [MAC]). Our collage, and the whole project of which it is
a part, is a brief
self--presentation of ethnomathematics of Western
culture, observed from the vantage point of 
the second half of the twentieth century.

\smallskip

Probably the most interesting intracultural 
interactions involving mathematics are those
that are not direct but rather proceed  via the mediation 
of value systems. A value system influences activities
in each domain and practically determines their cultural
interpretation. Conversely, an emerging value system
in one part of cultural activity (e.g. scientific) starts a process of reconsideration of other ones, their reformation,
sometimes leading to their extinction or total
remodeling. 

\smallskip

This is why in the last section I briefly touch upon
human values in the context of mathematical creativity.

\medskip

{\bf III.2. Rationality.}
Let us listen again to J.~R.~Newman (Introduction to vol. I of [WM]):

\smallskip

``... I began gathering the material for an anthology 
which I hoped would convey something of the diversity,
the utility and the beauty of mathematics''.

\smallskip

The book [WM] ``... presents mathematics as a tool, a language 
and a map; as a work of art and an end to itself; as a
fulfillment of the passion for perfection. It is seen as an object of
satire, a subject for humor and a source of controversy;
as a spur to wit and a leaven to the storyteller imagination; as an activity
which has driven men to frenzy and provided them with delight.
It appears in broad view as a body of knowledge
made by men, yet standing apart and independent of them.''

\smallskip

In this private and emotional list of values associated with mathematics
one is conspicuously absent: {\it rationality.}
One possible explanation is that in the Anglo--Saxon tradition, this
basic value of the Enlightenment came to be associated
with economic behavior, and often gets a narrow
interpretation: a rational actor is the one that 
consistently promotes self--interest.

\smallskip

Another explanation is that
being rational is not really delightful:
``Cogito ergo sum'' is an existence proof but it lacks
the urgency which a living soul feels without thinking.

\smallskip
 
Still, rationality in the Renaissance sense,
 ``Il natural desiderio di sapere'' (cf. [Ce]),
and the drive to be consistently rational
is a force without which the existence of mathematics 
through the centuries, and its successes in bringing its share
to the technological progress of society would be
impossible.

\medskip

{\bf III.3. Truth.} Extended and complex, subtle and mutually
contradictory views were expounded on the 
problem ``truth in mathematics'': see [Tr] for
a fairly recent review. Here I simply
state that axiologically, this is one of the central
values associated with mathematics, whatever
its historical and philosophical correlates might be.

\smallskip

Authority, practical efficiency, success
in competition, faith,
all these clashing values must recede in the mind of
a mathematician when he or she sets down to do their job.
  
\medskip

{\bf III.4. Action and contemplation.} By the nature of their trade,
mathematicians are inclined more to contemplation than to action.

\smallskip

The Romans, who were actors {\it par excellence} and revered Greek culture,
skipped Greek mathematics.
The imperial list of virtues -- valor, honor, glory, service --
did not leave much place for geometry.

\smallskip

This tradition continued through centuries, but as with any
tradition, there were exciting exceptions, and I will
conclude this  essay with a sketch of a great
mathematician of the last century, John von Neumann.

\smallskip

Neumann J\`anos was born on October 28, 1903 in Budapest,
and died in Washington, D.C., on February 8, 1957.
During this relatively short life span, he participated in,
and made crucial contributions to: the foundations of set theory,
quantum statistics and ergodic theory, game theory as a paradigm
of economic behavior, theory of operator algebras, 
the architecture of modern computers, the 
implosion principle for the creation of the hydrogen bomb,
and much more. 

\smallskip

Here are two samples of his thinking and modes of expression,
marking the beginning and the end of his career.

\smallskip

{\it Contemplation: The von Neumann Universe.} Cantor's description
of a set as an arbitrary collection of distinct elements of our  
thought is too generous in many contexts,
and the von Neumann Universe consists only of sets whose
elements are also sets. The potentially dangerous self--referentiality
is avoided by postulating that any family of sets $X_i$
such that $X_i$ is an element of $X_{i+1}$ has  a least element;
and the ultimate set, the least of all, is empty. 
Thus von Neumann Universe is
born from a ``philosophical vacuum'': its first elements are
$\emptyset$ (the empty set), $\{\emptyset\}$ (one--element set whose
only element is the empty set), $\{\{\emptyset\}\}$, $\{\emptyset , 
\{\emptyset\}\}$ etc.  Stingy curly brackets replace Cantor's
{\it Zusammenfassung \dots zu einem Ganzen}, and this operation,
which can be iteratively repeated, is the only one
that produces new sets from the already constructed ones.
Iteration can be, of course, transfinite, which was another
great insight of Cantor's.

\smallskip

It is difficult to imagine a purer object of contemplation than
this quiet and powerful hierarchy.

\smallskip

{\it Action: Hiroshima.} Excerpts from von Neumann's 
letter to R.~E.~Duncan,
IBM War History section, dated December 18, 1947 ([Neu], pp. 111--112):

\smallskip

` Dear Mr. Duncan 

In reply to your letter of December 16, [...] I can tell you the following things: I did initiate and carry out work during the war on
oblique shock reflection. This did lead to the conclusion
that large bombs are better detonated at a considerable altitude
than on the ground, since this leads to the higher
oblique--incidence pressure referred to. [...]

I did receive the Medal for Merit (October, 1946) and the
Distinguished Service Award (July, 1946). The citations
are as follows:

\smallskip

\centerline{``Citation to Accompany the Award of}

\centerline{The Medal for Merit}

\centerline{to}

\centerline{Dr. John von Neumann}

\smallskip

DR. JOHN VON NEUMANN, for exceptionally meritorious
conduct in the performance of outstanding services to the United
States from July 9, 1942 to August 31, 1945.
Dr. von Neumann, by his outstanding devotion to duty,
technical leadership, untiring cooperativeness, and sustained enthusiasm,
was primarily responsible for fundamental research by
the United States Navy on the effective use of high
explosives, which has resulted in the discovery of a new
ordnance principle for offensive action, and which has already
been proved to increase the efficiency of air power
in the atomic bomb attacks over Japan. His was a contribution
of inestimable value to the war
effort of the United States.

HARRY TRUMAN''

[...]'

\bigskip

\centerline{\bf References}

\medskip

[At] M.~Atiyah. {\it Geometry and Physics of the 20th Century.}
In: [GeoXX], pp. 4--9.

[B-BH] B.~Boo{\ss}--Bavnbek, J.~H{\o}yrup. {\it Introduction.}
In: [MW], pp. 1--19.

[Bour] N. Bourbaki. {\it Elements of the History of Mathematics.}
Springer Verlag, 1994.

[CaBa] M.~Li Calzi, A. Basile. {\it Economists and Mathematics
from 1494 to 1969. Beyond the Art of Accounting.}
In: [MC], pp. 95 -- 107.

[Ce] F.~Cesi. {\it Il natural desiderio di sapere.}
The Pontifical Academy of Sciences, Vatican, 2003.

[Cha] Yu.~V.~Chaikovski. {\it What is probability?
(Evolution of the notion from antiquity to Poisson)}
(Russian). In: Istoriko--Mathematicheskie Issledovaniya
(Studies in the History of Mathematics), ser. 2, 6(41),
Moscow 2001, pp. 34--56

[Che] K.~Chemla. {\it History of Mathematics in China: A Factor in World
History and a Source for New Questions.} In: Proceedings of the Int. Congr.
of Mathematicians, Berlin 1998, vol. III, pp. 789--798.
Documenta Mathematica, Bielefeld, 1998.

[CVS] {\it The Cultural Values of Science.} Proceedings of
the plenary session of the Pontifical Academy
of Sciences, Nov. 8 -- 11 2002. Vatican, 2003. 

[DaOl] H.~G.~Dales, G.~Oliveri. {\it Truth and the 
foundations of mathematics. An introduction.} In: [Tr], pp. 1--37.

[DH] P.~Davis, R.~Hersh. {\it The Mathematical Experience.} 
Birkh\"auser Boston, 1980.

[Ens] H.~M.~Enszenberger. {\it Drawbridge Up. Mathematics -- a Cultural Anathema.} A.~K.~Peters, Natick, Mass., 1999.

[GeoXX] {\it G\'eom\'etrie au XXe Si\`ecle. Histoire et
Horizons.} Ed. by  J.~Kouneiher, D.~Flament,
Ph.~Nabonnand, J.-J.~Szczeciniarz.)
Hermann, Paris, 2005. 


[KGSIPW] V.~Keilis-Borok, D.~Gascon, A.~Soloviev, M.~Intriligator,
R.Pichardo, F.~Winberg. {\it On the predictability
of crime waves in megacities -- extended summary.}
In: [CVS], pp. 221--229.

[KoMa] I.~Yu.~Kobzarev, Yu.~I.~Manin. {\it
Elementary Particles. Mathematics, Physics and Philosophy.}
Kluwer, 1989.


[MAC] {\it Mathematics Across Cultures. The history of
Non--Western Mathematics.} Ed. by H.~Selin. Kluwer
Academic Publishers, 2000.

[Macr] N.~Macrae. {\it John von Neumann.} AMS, 1999.

[Man1] Yu.~Manin.  {\it Truth, rigor and common sense.} 
In: [Tr], pp. 147--159.

[Man2]  Yu.~Manin.  {\it  Mathematics as metaphor.} In: Proc. of ICM, Kyoto 1990, vol. II, 1665--1671. The AMS and Springer Verlag.

[Man3]  Yu.~Manin.  {\it Georg Cantor and his heritage.} 
In: Algebraic Geometry: Methods, Relations, and Applications:
Collected papers dedicated to the memory
of Andrei Nikolaevich Tyurin. Proc. V. A. Steklov Inst. Math.
Moscow,  vol. 246, MAIK Nauka/Interperiodica, 2004,
195--203. Preprint math.AG/0209244

[MandHu] B.~Mandelbrot, R.~Hudson. {\it The (Mis)behaviour of Markets.}
Profile, 2005. 328 pp.

[MC] {\it Mathematics and Culture I.} Ed. by M. Emmer.
Springer 2004.


[Mu1] D.~Mumford. {\it The dawning of the age of stochasticity.}
In: Mathematics: Frontiers and Perspectives, ed. by V.~Arnold,
M.~Atiyah, P.~Lax and B.~Mazur. AMS, 2000, pp. 197--218.

[Mu2] D.~Mumford. {\it Pattern Theory: the Mathemaics of
Perception.} In: Proceedings of the Int. Congr.
of Mathematicians, Beijing 2002, vol. I, pp. 401--422.
Higher Education Press, Beijing, 2002.

[MW] {\it Mathematics and War.} Ed. by B.~Boo{\ss}--Bavnbek, J.~H\/oyrup.
Birkh\"auser Verlag, Basel, 2003.

[Neu] J.~von Neumann. {\it Selected Letters.}
ed. by M.~R\'edei. History of Math., vol 27,  AMS and London MS, 2005.

[New] J.~R.~Newman. {\it The Rhind Papyrus.}
In: [WM], vol. 1, pp. 170 -- 178.

[N{\o}] T.~N{\o}rretranders. {\it The User Illusion:
Cutting Conscioussness Down to Size.} Penguin, 1998.


[Pl] Plato. {\it Republic, Book VII.} 522d, 525c, 526 b--d .

[Po1] M.~Poovey. {\it A History of the Modern Fact.
Problem of Knowledge in the Sciences of Wealth and Society.}
The University of Chicago Press, 1998.

[Po2] M.~Poovey. {\it Can numbers ensure honesty?
Unrealistic expectations and the US accounting scandal.}
Notices of the AMS, vol. 50:1, Jan. 2003, pp. 27--35.

[Qu] Anjing Qu. {\it The third approach to the history
of mathematics in China.} In: Proceedings of the Int. Congr.
of Mathematicians, Beijing 2002, vol. III, pp. 947--958.
Higher Education Press, Beijing, 2002.




[Sci] {\it Bourbaki. Une soci\'et\'e secr\`ete de math\'ematiciens.}
Pour la Science, No 2, 2000.

[S-S] R.~Siegmund--Schultze. {\it Military Work in Mathematics 1914--1945: an Attempt at an International Perspective.} In: [MW], pp. 23--82.

[Sol] R.~Solow. {\it How did economics get that way
and what way did it get?} Daedalus, Fall 2005, pp. 87--100.

[So] S.~Sontag. {\it Illness as Metaphor, and AIDS and its Metaphors.}
Picador; Farrar, Strauss and Giroux, NY 1990.

[Tr] {\it Truth in Mathematics.} Ed. by H.~G.~Dales and G.~Oliveri.
Clarendon Press, Oxford, 1998.

[Tu] H.~W.~Turnbull. {\it The great mathematicians.}
In: [WM], vol. 1, pp. 75--168.

[WM] {\it The World of Mathematics.} A small library of the
literature of mathematics from A'h--mos\'e the Scribe to
Albert Einstein, presented with commentaries and
notes by James R. Newman. Vols 1--3.  Simon and Schuster,
New York, 1956.

\enddocument